# New Lower Bounds for the Number of Blocks in Balanced Incomplete Block Designs

MUHAMMAD A. KHAN

ABSTRACT: Bose [1] proved the inequality $b \geq v + r - 1$ for resolvable balanced incomplete block designs (RBIBDs) and Kageyama [4] improved it for RBIBDs which are not affine resolvable. In this note we prove a new lower bound on the number of blocks $b$ that holds for all BIBDs. We further prove that for a significantly large number of BIBDs our bound is tighter than the bounds given by the inequalities of Bose and Kageyama.

## 1. INTRODUCTION

Let $D$ be a balanced incomplete block design (BIBD) with parameters $(v, b, r, k, \lambda)$. The following results follow immediately from definition.

$$bk = v\lambda, \qquad r(k-1) = \lambda(v-1)$$

A less obvious relation is Fisher's inequality [3].

$$b \geq v \quad \text{(or equivalently } r \geq k\text{)}$$

In his classical paper on combinatorial designs Bose [1] strengthened Fisher's inequality for RBIBDs proving that

$$b \geq v + r - 1$$

Later Michael [6] proved the same inequality under the weaker assumption $v = nk$, for some positive integer $n$. Stanton [7] showed that a necessary and sufficient condition for Bose's inequality is $r \geq \lambda + k$ (which is weaker than $v = nk$).

Over the years Bose's inequality has been improved under stronger assumptions. For instance Kageyama [4] showed that $b \geq 2v + r - 1$ for resolvable designs which are not affine resolvable. However, so far the lower bound given by Bose's inequality has not been improved generally. The question that motivated this paper is the following:

*"Can we establish a lower bound for b that works for all BIBDs and at the same time improves Bose's bound for a significant number of BIBDs?"*

Such results are important as they can help in answering questions regarding the existence of BIBDs. Here we prove the following theorems:

---

**Theorem 1.1** *For a BIBD with parameters* $(v, b, r, k, \lambda)$

a) *if D is nontrivial then* $b \geq \left\lceil \dfrac{(v-k)^3}{v^2} \right\rceil + 2r - \lambda$.

b) *if D is nontrivial and* $r \geq v - 1$ *then* $b \geq \left\lceil \dfrac{(v-k)^2}{v-1} \right\rceil + 2r - \lambda$.

We prove theorem 1.1 in section 2, establishing some new inequalities for BIBD s along the way. In section 3 we compare the lower bounds given in theorem 1.1 with Bose's lower bound for fifty admissible parameter sets listed in the Handbook of Combinatorial Designs [2]. We also demonstrate that even in the special case of resolvable designs which are not affine resolvable our bounds are tighter than Kageyama's bound.

## 2. PROOF OF THEOREM 1.1 AND RELATED RESULTS

A BIBD with $k = 1$ or $k = v$ is trivial. So $1 < k < v$ for any nontrivial BIBD. We begin with some preliminary but significant results.

**Lemma 2.1** *Let D be a BIBD with parameters* $(v, b, r, k, \lambda)$. *Then*
  a) $bk^2 \geq \lambda v^2$.
  b) *if D is nontrivial then* $bk^m < \lambda v^m$ *for* $m > 3$. *In particular* $bk^3 < \lambda v^3$.
  c) *if D is nontrivial then* $k^m < \lambda v^{m-1}$ *for* $m > 3$. *In particular* $k^3 < \lambda v^2$.

*Proof.*

a) $bk^2 - \lambda v^2 = \dfrac{vr}{k} k^2 - \lambda v^2 = vk \dfrac{\lambda(v-1)}{k-1} - \lambda v^2$

$= \dfrac{\lambda v^2 k - \lambda vk - \lambda v^2 k + \lambda v^2}{k-1} = \dfrac{\lambda v^2 - \lambda vk}{k-1} \geq 0$

b) It suffices to prove the result for $m = 3$. Observe that the function $f(x) = \dfrac{x^2}{x-1}$ is strictly increasing for $x \geq 2$. So for $k \geq 2$

$\dfrac{k^2}{k-1} < \dfrac{v^2}{v-1} \quad \Rightarrow \quad \dfrac{\lambda v(v-1)}{k(k-1)} k^3 < \lambda v^3$

$\quad \Rightarrow \quad bk^3 < \lambda v^3$

But $D$ is nontrivial so we always have $k \geq 2$.

c) Follows from b) and Fisher's inequality.

A similar argument leads to the corresponding result for 1-designs.

**Lemma 2.2** *If $D$ is a 1-$(V, B, R, K, \Lambda)$ design then*
  a) *for $m > 2$ we have $BK^m \leq \Lambda V^m$. In particular $BK^2 \leq \Lambda V^2$.*
  b) *for $m > 2$ and $B \geq V$ we have $K^m \leq \Lambda V^{m-1}$. In particular $K^2 \leq \Lambda V$.*

*Proof.* Follows from $BK = RV = \Lambda V$ and $V \geq K$.

Note that for lemma 2.2 b) we explicitly need to assume $b \geq v$ as this does not hold in general for 1-designs.

The question may be asked as to whether $k^2 \leq \lambda v$ holds for a BIBD with parameters $(v, b, r, k, \lambda)$. Clearly the bound in theorem 1.1 can be improved if we replace $k^3 \leq \lambda v^2$ with $k^2 \leq \lambda v$. In fact a direct proof of $k^2 \leq \lambda v$ together with lemma 2.1 b) would prove Fisher's inequality. However neither $k^2 \leq \lambda v$ nor $k^2 \geq \lambda v$ is obvious from lemma 2.1. Here we shall use a variant of lemma 2.2 to establish the relation $k^2 \leq \lambda v$ for a large class of BIBDs. The proof uses the idea of leave of a $t$-design. Let $E$ be a $t$-$(v, k, \lambda)$ design with $b$ blocks and replication number $r$. The leave $E^L$ of $E$ is a $(t-1)$-$(v-1, k, \lambda)$ design with $(b-r)$ blocks (see [2, p 120] for complete definition).

**Theorem 2.3** *If $D$ is a BIBD with parameters $(v, b, r, k, \lambda)$ and $b \geq v + r - 1$ then $k^2 \leq \lambda(v-1)$.*

*Proof.* The leave $D^L$ of $D$ is a 1-design. Since $b - r \geq v - 1$ therefore we can apply lemma 2.2 b) to $D^L$ so that $k^2 \leq \lambda(v-1)$.

The above proof provides an instance where 1-designs have been used to prove a theorem about BIBDs (2-designs). This is extremely rare in literature as BIBDs have a rich structure compared to 1-designs.

Now we can return to our main results. For the proof of theorem 1.1 recall that the complementary design of a BIBD with parameters $(v, b, r, k, \lambda)$ is a BIBD with parameters $(v, b, b - r, v - k, b - 2r + \lambda)$ provided $b - 2r + \lambda \neq 0$.

*Proof of Theorem 1.1.*
a) Let $D$ be a BIBD satisfying the hypothesis. As $D$ is nontrivial therefore $v \neq k$ and

$$b - 2r + \lambda = \frac{\lambda v(v-1)}{k(k-1)} - 2\frac{\lambda(v-1)}{(k-1)} + \lambda = \lambda\left[\frac{(v-1)}{(k-1)}\left(\frac{v}{k} - 2\right) + 1\right] \neq 0$$

So $D^C$ exists and has parameters $(v, b, b - r, v - k, b - 2r + \lambda)$. Now applying the case $m = 3$ of lemma 2.1 c) gives $(v-k)^3 < (b - 2r + \lambda)v^2$ and so $b > \frac{(v-k)^3}{v^2} + 2r - \lambda$. The result follows since $b$ is a positive integer.

b) Let $(v^*, b^*, r^*, k^*, \lambda^*) = (v, b, b - r, v - k, b - 2r + \lambda)$.

If $r \geq v - 1$ then $b - (b - r) \geq v - 1$ giving $b^* - r^* \geq v^* - 1$. Thus theorem 2.3 is applicable and we have $(v-k)^2 \leq (b - 2r + \lambda)(v-1)$ or $b \geq \dfrac{(v-k)^2}{v-1} + 2r - \lambda$.

## 3. SHARPNESS OF LOWER BOUNDS

The following table lists Bose's lower bound and the bounds given by theorem 1.1 for admissible parameter sets listed in the [2, pp 36-58] ordered lexicographically by $r$, $k$ and $\lambda$. The first column gives the serial number assigned in The Handbook of Combinatorial Designs. The entries in bold indicate that the bound(s) given by theorem 1.1 are more stringent than the one given by Bose's inequality. The dash specifies that the bound is not applicable.

**TABLE 1 Comparison of Lower Bounds**

| No | Admissible Parameter Sets | Theorem 1.1 Part (a) | Theorem 1.1 Part (b) | Bose $v + r - 1$ |
|---|---|---|---|---|
| 1 | (7,7,3,3,1) | 7 | - | - |
| 2 | (9,12,4,3,1) | 10 | - | 12 |
| 3 | (13,13,4,4,1) | 12 | - | - |
| 4 | (6,10,5,3,2) | 9 | **10** | 10 |
| 5 | (16,20,5,4,1) | 16 | - | 20 |
| 101 | (8,28,14,4,6) | **23** | **25** | 21 |
| 102 | (15,42,14,5,4) | **29** | **32** | 28 |
| 103 | (36,84,14,6,2) | **47** | - | 49 |
| 104 | (15,35,14,6,5) | 27 | **29** | 28 |
| 105 | (85,170,14,7,1) | 93 | - | 98 |
| 201 | (28,72,18,7,4) | 44 | - | 45 |
| 202 | (64,144,18,8,2) | 77 | - | 81 |
| 203 | (145,290,18,9,1) | 155 | - | 162 |
| 204 | (73,146,18,9,2) | 84 | - | 90 |
| 205 | (49,98,18,9,3) | 60 | - | 66 |
| 301 | (85,105,21,17,4) | 82 | - | 105 |
| 302 | (120,140,21,18,3) | 113 | - | 140 |

| No | Admissible Parameter Sets | Theorem 1.1 | | Bose |
| --- | --- | --- | --- | --- |
| | | Part (a) | Part (b) | $v + r - 1$ |
| 303 | (190,210,21,19,2) | 179 | - | 210 |
| 304 | (400,420,21,20,1) | 384 | - | 420 |
| 305 | (421,421,21,21,1) | 403 | - | - |
| 401 | (529,552,24,23,1) | 510 | - | 552 |
| 402 | (553,553,24,24,1) | 532 | - | - |
| 403 | (277,277,24,24,2) | 258 | - | - |
| 404 | (185,185,24,24,3) | 167 | - | - |
| 405 | (139,139,24,24,4) | 123 | - | - |
| 501 | (55,99,27,15,7) | 69 | - | 81 |
| 502 | (460,690,27,18,1) | 462 | - | 486 |
| 503 | (153,231,27,18,3) | 158 | - | 180 |
| 504 | (52,78,27,18,9) | 60 | - | 78 |
| 505 | (91,117,27,21,6) | 90 | - | 117 |
| 601 | (61,366,30,5,2) | **106** | - | 90 |
| 602 | (41,246,30,5,3) | **85** | - | 70 |
| 603 | (31,186,30,5,4) | **75** | **79** | 60 |
| 604 | (25,150,30,5,5) | **68** | **72** | 54 |
| 605 | (21,126,30,5,6) | **64** | **67** | 50 |
| 701 | (17,68,32,8,14) | **53** | **56** | 48 |
| 702 | (145,464,32,10,2) | **180** | - | 176 |
| 703 | (25,80,32,10,12) | **58** | **62** | 56 |
| 704 | (33,96,32,11,10) | **64** | **70** | 64 |
| 705 | (177,472,32,12,2) | 206 | - | 208 |
| 801 | (69,138,34,17,16) | 90 | - | 102 |
| 802 | (35,70,34,17,16) | 57 | 62 | 68 |
| 803 | (715,1105,34,22,1) | 719 | - | 748 |
| 804 | (69,102,34,23,11) | 78 | - | 102 |

| No | Admissible Parameter Sets | Theorem 1.1 | | Bose |
| --- | --- | --- | --- | --- |
| | | Part (a) | Part (b) | $v + r - 1$ |
| 805 | (154,187,34,28,6) | 147 | - | 187 |
| 901 | (13,78,36,6,15) | **60** | **62** | 48 |
| 902 | (217,1116,36,7,1) | **268** | - | 252 |
| 903 | (28,144,36,7,8) | **76** | **81** | 63 |
| 904 | (64,288,36,8,4) | **111** | - | 99 |
| 905 | (22,99,36,8,12) | **66** | **70** | 57 |

Looking at table 1 it seems that very few BIBDs satisfy the condition $r \geq v - 1$. But the following theorem shows that this is not the case. Indeed if a nontrivial BIBD satisfies Bose's inequality then either $D$ or $D^C$ must satisfy this condition.

**Theorem 3.1** *If $D$ is a nontrivial BIBD with parameters $(v, b, r, k, \lambda)$ such that $b \geq v + r - 1$ and $r < v - 1$. Suppose the complementary design $D^C$ has parameters $(v^*, b^*, r^*, k^*, \lambda^*)$ then $r^* \geq v^* - 1$.*

*Proof.* Here $(v^*, b^*, r^*, k^*, \lambda^*) = (v, b, b - r, v - k, b - 2r + \lambda)$. Now $b - r \geq v - 1$ so that $r^* \geq v^* - 1$.

It is worth noting that the lower bound in part a) of theorem 1.1 is satisfied by all nontrivial BIBDs regardless of their parameters. This is quite unique as only Fisher's inequality is known to offer such a general meaningful lower bound. At the same time both of our inequalities improve Bose's inequality for several BIBDs. In particular as $r$, $k$ and $\lambda$ increase lexicographically our bounds become considerably better than Bose.

Interestingly even Kageyama's bound is not always stronger than any of our lower bounds. For instance consider the RBIBD with parameters (16, 140, 35, 4, 7) which is not affine resolvable [4, Example (iii)]. The inequalities in theorem 1.1 a) and b) respectively say $b \geq 71$ and $b \geq 73$ whereas Kageyama's inequality only says $b \geq 65$.

**Concluding Remarks:** In this paper we have established two new lower bounds for the number of blocks in a balanced incomplete block design which improve all known lower bounds for a large number of BIBDs. At the same time we have suggested a possible approach for proving Fisher's inequality combinatorially without using linear algebra. Currently we are extending the methods described here to the more general setting of *t*-designs [5]. There is plenty of room for improvement even for the lower bounds proved for BIBDs. We conjecture the following:

**Open Problem 3.1** *If D is a nontrivial BIBD with parameters (v, b, r, k, λ) and $b \geq v + r - 1$ then $b \geq \left\lceil \frac{(v-k)^2}{v-1} \right\rceil + 2r - \lambda$.*

So far the computational evidence supports this conjecture.

**Acknowledgements:** The author is thankful to King Fahd University of Petroleum & Minerals for continuous research support.

PREPARATORY YEAR MATHEMATICS PROGRAM
KING FAHD UNIVERSITY OF PETROLEUM & MINERALS
DHAHRAN 31261, SAUDI ARABIA

Email: malikhan@kfupm.edu.sa